\documentstyle[art12]{article}
\title{Weierstrass Gap Sequence at Total Inflection Points of Nodal Plane
 Curves}
\author{Marc Coppens \and Takao Kato}
\let\bbb=\bf
\def\al{\alpha}
\def\be{\beta}
\def\ga{\gamma}
\def\Ga{{\mit \Gamma}}
\def\de{\delta}
\def\fracd#1#2{{\displaystyle{#1}\over\displaystyle{#2}}}
\def\sus{\subset}
\def\vs#1{\par\vspace{.#1in}}
\date{}
\begin{document}
\newtheorem{lem}{Lemma}[section]
\newtheorem{thm}[lem]{Theorem}
\newtheorem{prop}[lem]{Proposition}
\newtheorem{cor}[lem]{Corollary}
\maketitle
\setcounter{section}{-1}

\section{Introduction}

Let $\Ga$ be a plane curve of degree $d$ with $\de$ ordinary nodes and no
other singularities. Let $C$ be the normalization of $\Ga$. Let
$g=\fracd{(d-1)(d-2)}{2}-\de$; the genus of $C$. We identify smooth
points of $\Ga$ with the corresponding points on $C$. In particular, if
$P$ is a smooth point on $\Ga$ then the Weierstrass gap sequence at $P$
 is considered with respect to $C$. A smooth point $P\in\Ga$ is called an
 $(e-2)$-inflection point if $i(\Ga ,T;P)=e\ge 3$ where $T$ is the
tangent line to $\Ga$ at $P$
(cf. Brieskorn--Kn\"{o}rrer\cite[p.~372]{B}). Of course, $e\le d$ and a
1-inflection point is an ordinary flex. In particular, a
$(d-2)$-inflection point is called a total inflection point.

Let $N$ be the semigroup consisting of the non-gaps of $P$, so
${\bbb N}-N=\{\al_1<\al_2<\cdots <\al_g\}$ is the Weierstrass gap
sequence of $P$. Clearly $\{ d-1,d\}\sus N$, so
$N_d:=\{ a(d-1)+bd|a,b\in{\bbb N}\}\sus N$ (see also Lemma \ref{lem:1}).

Let $k=\min\{\ell\in{\bbb N}|\de\le\fracd{\ell (\ell +3)}{2}\}$ and let
$$
N_{d,\de}^{(1)}=N_d\cup\{ n\in{\bbb N}|n\ge (d-k-3)d+\fracd{k(k+3)}{2}
-\de +2\}.
$$
Let ${\bbb N}-N_{d,\de}^{(1)}=\{\al_1^{(1)}<\al_2^{(1)}<\cdots <
\al_g^{(1)}\}$. One has $\al_i\ge\al_i^{(1)}$ for $1\le i\le g$. So
$N_{d,\de}^{(1)}$ is the minimal (in the sense of weight) possible
semigroup of non-gaps.

For $\de\in\{ 0,1\}$, one has $N=N_{d,\de}^{(1)}$. For $\de\ge 2$ there
exist pairs of $(\Ga ;P)$ as above with $N\ne N_{d,\de}^{(1)}$.
We give a list of all possible values for $N$ in case $2\le\de\le 5$.
(see end of \S 1).

Define $N_{d,1}^{({\rm max})}=N_{d,1}^{({\rm max2})}=N_{d,1}$ and, by
means of induction, for $\de\ge 2$,
\begin{eqnarray*}
N_{d,\de}^{({\rm max})} & = & N_{d,\de -1}^{({\rm max})}\cup
\{ (d-\de -2)d+1\} \\
N_{d,\de}^{({\rm max2})} & = & N_{d,\de -1}^{({\rm max2})}\cup
\{ (d-\de -2)d+\de\}.
\end{eqnarray*}
$N_{d,\de}^{({\rm max})}$ (resp. $N_{d,\de}^{({\rm max2})}$) is a
semigroup if and only if $d\ge 2\de +1$ (resp. $2\de$). Let
\begin{eqnarray*}
{\bbb N}-N_{d,\de}^{({\rm max})} & = &
\{\al_1^{({\rm max})}<\al_2^{({\rm max})}<\cdots <\al_g^{({\rm max})}\} \\
{\bbb N}-N_{d,\de}^{({\rm max2})} & = &
\{\al_1^{({\rm max2})}<\al_2^{({\rm max2})}<\cdots <\al_g^{({\rm max2})}\}.
\end{eqnarray*}
We prove that $\al_i\le\al_i^{({\rm max})}$ for $1\le i\le g$ and if
$N\ne N_{d,\de}^{({\rm max})}$, then $\al_i\le\al_i^{({\rm max2})}$ for
$1\le i\le g$ (Lemma \ref{lem:max}). So $N_{d,\de}^{({\rm max})}$ (resp.
$N_{d,\de}^{({\rm max2})}$) is the maximal (resp. up to 1 maximal)
semigroup of non-gaps.

Our main results are the following:
\begin{enumerate}
\item   There exist pairs $(\Ga ;P)$ such that $N=N_{d,\de}^{(1)}$
(\ref{thm:general}),

\item   If $d\ge 2\de +1$ (resp. $d\ge 2\de$) then there exist pairs
$(\Ga ;P)$ such that $N=N_{d,\de}^{(max)}$ (resp. $N=N_{d,\de}^{(max2)}$)
 (\ref{prop:max}).
\end{enumerate}

The existence of Weierstrass points with gap sequence
${\bbb N}-N_{d,\de}^{(1)}$ is already proved in \cite{K} for the case
$\de =\fracd{d^2-7d+12}{2}$. The method used in that paper is completely
different from ours. It has the advantage of not using plane models but
the proof looks more complicated. It might be possible to prove our
existence result in this way completely, but it might become very
complicated. We didn't try it. Also, it gives an affirmative answer to
Question 1 in \cite{C} for the case $s=n+1$. It is not clear to us at the
 moment how to generalize the proof for the cases with $s\ge n+2$.

\section{Generalities and low values for $\de$}

To start, we deal with the case $\de =0$.
\vs 2
\begin{lem}
Let $\Ga$ be a smooth plane curve of degree $d$ and let $P$ be a
total inflection point of $\Ga$. Then $N_d=N_{d,0}^{(1)}$ is the
semigroup of non-gaps of $P$.
\label{lem:0}
\end{lem}
\vs 1
{\it Proof\/}. Let $T$ be the tangent line at $P$, $L_1$ be a general
line passing through $P$ and let $L_2$ be a general line not passing
through $P$. Then the curve $C(a,b)=aT+bL_1+(d-3-a-b)L_2$ is canonical
adjoint, if $0\le a\le d-3, 0\le b\le d-3-a$. Then we have
$i(\Ga .C(a,b);P)=ad+b$. Hence,
$\{ ad+b+1:0\le a\le d-3, 0\le b\le d-3-a\}$ is the gap sequence at $P$.
This completes the proof.
\vs 2
In order to study the case $\de >0$, we prove some lemmas. For the rest
of this section, $\Ga$ is a plane curve of degree $d$ with $\de (>0)$
ordinary nodes $s_1,\dots ,s_{\de}$ as its only singularities. Also
$P\in\Ga$ is a total inflection point.
\vs 2
\begin{lem}
The set of nongaps at $P$ contains $N_{d,0}$.
\label{lem:1}
\end{lem}
\vs 2
{\it Proof\/}. Assume that $n\in N_{d,0}$. Let $\al =\left[\fracd{n-1}{d}
\right] +1$, $\ell$ be the equation of $T$ (the tangent line at $P$),
$\ell_0$ be the equation of a general line passing through $P$ and let
$\ell_1$ be the equation of a general line. Considering
$$
\frac{\ell_0^{\al d-n}\ell_1^{\al +n-\al d}}{\ell^{\al}},
$$
we obtain that $n$ is a nongap at $P$.
\vs 2
%Next, we recall a lemma due to Namba \cite[Lemma 2.3.2]{namba}
% (cf. Coppens and Kato \cite[Lemma 1.1]{ck1} for a generalization).
%\vs 1
%\begin{lem}
%Let $C$ be a nonsigular plane curve defined by $f(x,y)=0$. Let $g(x,y),
%h(x,y)$ be polynomials such that $f, g, h$ have no common factor each
%other. Assume  $f(0,0)=0$. Then,
%$i(g,h;O)\ge\min\{ i(f,g;O), i(f,h;O)\}$.
%\label{lem:namba}
%\end{lem}
%\vs 1
%Using this lemma, we have
\begin{lem}
Let $\ga$ be a curve of degree less than $d$ so that
$i(\ga ,\Ga ;P)=k\ge d$. Then, $T$ is a component of $\ga$, i.~e. there
is a curve $\ga'$ of degree $\deg\ga -1$ such that $\ga =\ga'T$.
\label{lem:2}
\end{lem}
\vs 2
{\it Proof\/}. Since $i(T,\Ga ;P)=d$ and $i(\ga ,\Ga ;P)\ge d$, by
Namba's lemma \cite[Lemma 2.3.2]{namba} (cf. Coppens and Kato
\cite[Lemma 1.1]{ck1} for a generalization), we have
$i(T,\ga ;P)\ge d>\deg\ga$.
Hence we have the desired result by Bezout's theorem.
\vs 2
By a successive use of this lemma we have:
\vs 2
\begin{lem}
Let $\ga$ be a canonical adjoint curve such that
$i(\ga ,\Ga ;P)=\al d+\be$ $(0\le\al\le d-3, 0\le\be\le d-3-\al )$.
Then, there is an adjoint curve $\ga'$ of degree $d-3-\al$ such that
$\ga =T^{\al}\ga'$ and $i(\ga',\Ga ;P)=\be$.
\label{lem:3}
\end{lem}
\vs 2
Using Lemma \ref{lem:3} we have the following corollaries:
\vs 2
\begin{cor}
If $\de\ge 1$, then $i(\ga ,\Ga ;P)<(d-3)d$ for every canonical adjoint
curve $\ga$, hence $(d-3)d+1$ is a nongap at $P$.
\label{cor:1}
\end{cor}
\vs 2
\begin{cor}
Assume that $\de\ge 2$. Then, $(d-4)d+\be +1$ $(\be =0$ or $1)$ is a gap
if and only if there is a line $L_0$ such that
$s_1,\dots ,s_{\de}\in L_0$. Moreover, in this case, the following three
conditions are equivalent:
\begin{enumerate}
        \item   $P\notin L_0$, {\rm (resp.} $P\in L_0)$,

        \item   $(d-4)d+1$ {\rm (resp.} $(d-4)d+2)$ is a gap,

        \item   $(d-3-\al )d+1+\al$ $(\al =1,\dots ,\de -1)$ {\rm (resp.}
                $(d-3-\al )d+1$ $(\al =1,\dots ,\de -1))$ are nongaps.
\end{enumerate}
\label{cor:2}
\end{cor}
\vs 2
{\it Proof\/}. The existence of the line $L_0$ and the equivalence
between (i) and (ii) follows immediately from Lemma \ref{lem:3}.

Assume that $(d-4)d+1$ is a gap. If $(d-3-\al )d+\al +1$
$(1\le\al\le\de -1)$ is a gap then Lemma \ref{lem:3} provides an adjoint
curve $\ga'$ of degree $\al$ with $i(\ga' ,\Ga ;P)=\al$. So $\ga'$ has
$s_1,\dots ,s_{\de}$ as common points with $L_0$. Bezout's theorem
implies that $\ga'=\ga''L_0$ where $\ga''$ is a curve of degree $\al -1$
with $i(\ga'',\Ga ;P)=\al$ (since $P\notin L_0$). Namba's lemma implies
$\ga''=\ga'''T$, but then $i(\ga'',\Ga ;P)\ge d$, so
$\de\ge\al +1\ge d+1$. A contradiction since $s_1,\dots ,s_{\de}$ are
collinear.

Assume that $(d-4)d+2$ is a gap. If $(d-3-\al )d+1$ $(1\le\al\le\de -1)$
is a gap then Lemma \ref{lem:3} provides an adjoint curve $\ga'$ of
degree $\al$ with $i(\ga' ,\Ga ;P)=0$. But $\ga'=\ga''L_0$ and
$P\in L_0$, hence a contradiction.

Assuming (iii), we obtain (ii) because the number of gaps has to be $g$.
\vs 2
Using Lemma \ref{lem:1} and Corollary \ref{cor:2}, we are able to
determine the gap sequence in case that $s_1,\dots ,s_{\de}$ are
collinear.

Checking case by case by use of Lemmas \ref{lem:1} and \ref{lem:3}, we
show a table of possible nongaps $N_{d,\de}$ for $1\le\de\le 5$.

\newpage
{
\parindent -4.5em
\begin{tabular}[t]{|p{13em}|p{22.5em}|}
\hline
$N_{d,1}=N_{d,0}\cup\{ (d-3)d+1\}$ &
general\\
\hline\hline
$N_{d,2}^{(1)}=N_{d,1}\cup\{ (d-4)d+2\}$ & general\\
\hline
$N_{d,2}^{(2)}=N_{d,1}\cup\{ (d-4)d+1\}$ & $s_1, s_2, P$ are collinear\\
\hline\hline
$N_{d,3}^{(1)}=N_{d,2}^{(1)}\cup\{ (d-4)d+1\}$ & general\\
\hline
$N_{d,3}^{(2)}=N_{d,2}^{(1)}\cup\{ (d-5)d+3\}$ & $s_1, s_2, s_3$ are
collinear but not $P$\\
\hline
$N_{d,3}^{(3)}=N_{d,2}^{(2)}\cup\{ (d-5)d+1\}$ & $s_1, s_2, s_3, P$ are
collinear\\
\hline\hline
$N_{d,4}^{(1)}=N_{d,3}^{(1)}\cup\{ (d-5)d+3\}$ & general\\
\hline
$N_{d,4}^{(2)}=N_{d,3}^{(1)}\cup\{ (d-5)d+2\}$ & $s_1,\dots ,s_4$ general
 but $i(\ga ,\Ga ;P)=2$ where $\ga$ is the conic passing through
$s_1,\dots ,s_4,P$\\
\hline
$N_{d,4}^{(3)}=N_{d,3}^{(1)}\cup\{ (d-5)d+1\}$ & $s_1,s_2,s_3,P$ are
collinear but not $s_4$\\
\hline
$N_{d,4}^{(4)}=N_{d,3}^{(2)}\cup\{ (d-6)d+4\}$ & $s_1,s_2,s_3,s_4$ are
collinear but not $P$\\
\hline
$N_{d,4}^{(5)}=N_{d,3}^{(3)}\cup\{ (d-6)d+1\}$ & $s_1,s_2,s_3,s_4,P$ are
collinear\\
\hline\hline
$N_{d,5}^{(1)}=N_{d,4}^{(1)}\cup\{ (d-5)d+2\}$ & general\\
\hline
$N_{d,5}^{(2)}=N_{d,4}^{(1)}\cup\{ (d-5)d+1\}$ & $s_1,\dots ,s_5$ general
 but $\exists$ conic $\ga$ passing through $s_1,\dots ,s_5, P$ and
$i(\ga ,\Ga ;P)=1$\\
\hline
$N_{d,5}^{(3)}=N_{d,4}^{(2)}\cup\{ (d-5)d+1\}$ & $s_1,\dots ,s_5$ general
 but $\exists$ conic $\ga$ passing through $s_1,\dots ,s_5, P$ and
$i(\ga ,\Ga ;P)=2$\\
\hline
$N_{d,5}^{(4)}=N_{d,4}^{(1)}\cup\{ (d-6)d+4\}$ & $s_1,\dots ,s_4$ are
collinear but not $s_5, P$\\
\hline
$N_{d,5}^{(5)}=N_{d,4}^{(3)}\cup\{ (d-6)d+1\}$ & $s_1,\dots ,s_4, P$ are
collinear but not $s_5$\\
\hline
$N_{d,5}^{(6)}=N_{d,4}^{(4)}\cup\{ (d-7)d+5\}$ & $s_1,\dots ,s_5$ are
collinear but not $P$\\
\hline
$N_{d,5}^{(7)}=N_{d,4}^{(5)}\cup\{ (d-7)d+1\}$ & $s_1,\dots ,s_5, P$ are
collinear\\
\hline\hline
\end{tabular}
}

\newpage
\section{General Case ($\de\ge 2$)}
Remember the definition of $N_{d,\de}^{(1)}$, let
$k=\min\{\ell\in{\bbb N}|\de\le\fracd{\ell (\ell +3)}{2}\}$.
Then
$$
N_{d,\de}^{(1)}=N_d\cup\{ n\in{\bbb N}|n\ge (d-k-3)d+\fracd{k(k+3)}{2}
-\de +2\}.
$$

In this section, we prove that for $(\Ga ;P)$ general, the semigroup of
non-gaps of $P$ is equal to $N_{d,\de}^{(1)}$.
\vs 1
Let ${\bbb P}_{\ell}\cong{\bbb P}^{\ell (\ell +3)/2}$ be the linear system of
divisors of degree $\ell$ on ${\bbb P}^2$. Let
$$
{\bbb P}_{\ell}(s_1,\dots ,s_{\de})=\{\ga\in{\bbb P}_{\ell}|s_1,\dots ,
s_{\de}\in\ga\},
$$
and let
$$
{\bbb P}_k(s_1,\dots ,s_{\de};m)=\{\ga\in{\bbb P}_k(s_1,\dots ,
s_{\de})|i(\Ga ,\ga ;P)\ge m\}.
$$
\vs 1
\begin{lem}
Assume that
$$
(*)\qquad\left\{
\begin{array}{ll}
{\bbb P}_{\ell}(s_1,\dots ,s_{\de})=\emptyset\quad & {\rm if\ }\ell <k,\\
{\bbb P}_k(s_1,\dots ,s_{\de};m)=\emptyset\quad & {\rm if\ }m>
\fracd{k(k+3)}{2}-\de .
\end{array}\right.
$$
Then the Weierstrass gap sequence of $\Ga$ at $P$ is given by
${\bbb N}^+-N^{(1)}_{d,\de}$.
\label{lem:general}
\end{lem}
\vs 2
{\it Proof\/}. By Lemma \ref{lem:1}, every element of $N_{d,0}$ is a
nongap. For $0\le n\le d-3$ the natural number not belonging to
$N_{d,0}$ are $nd+1,\dots ,nd+(d-n-2)$. Assume such a number $nd+\be$
(hence $0\le n\le d-3$, $1\le\be\le d-n-2$) is a gap. Then there exists a
 canonical adjoint curve $\ga$ of $\Ga$ with
$$
i(\ga ,\Ga ;P)=nd+\be -1.
$$
Lemma \ref{lem:3} gives us that there exists
$\ga'\in{\bbb P}_{d-3-n}(s_1,\dots ,s_{\de})$ with
$i(\ga',\Ga ;P)=\be -1$. But the hypothesis $(*)$ implies that this is
impossible for $d-3-n<k$, i.e. $n>d-3-k$ or for $n=d-3-k$ and
$\be -1>\fracd{k(k+3)}{2}-\de$. So, the only possible gaps are
$$
\begin{array}{llcl}
1, \ \ \ 2, & \dots & \dots ,& d-2\\
d+1, \ d+2, & \dots & \dots ,& 2d-3\\
2d+1, 3d+2, & \dots & \dots ,& 3d-4\\
 & \dots & \dots & \\
(d-4-k)d+1, & \dots & \dots ,& (d-3-k)d-(d-2-k)\\
(d-3-k)d+1, & \dots & \dots ,& (d-3-k)d+\fracd{k(k+3)}{2}-\de +1.
\end{array}
$$
Since these are $g$ mumbers, we obtain the gaps of $C$ at $P$. It is
clear that this set is ${\bbb N}^+-N^{(1)}_{d,\de}$.
\vs 2
\begin{thm}
The hypothesis $(*)$ in Lemma {\rm \ref{lem:general}} occurs.
\label{thm:general}
\end{thm}
\vs 1
{\it Proof\/}.
(Inspired by the proof of Proposition 3.1 in \cite{T}).
Take a union of $d$ general lines in ${\bbb P}^2$: $\Ga_0=L_1\cup L_2\cup
\cdots\cup L_d$.

Let $P_1=L_1\cap L_2$, $\{ P_2,P_3\} =L_3\cap (L_1\cup L_2)$ and so on.
Take $0\le\de\le\fracd{(d-1)(d-2)}{2}$. The statement $(*)$ holds for
$\Ga_0$ instead of $\Ga$ and $s_1=P_1,\dots ,s_{\de}=P_{\de}$ and $P_0$
suitably chosen on $L_d$.

Indeed, let $k=\min\{\ell\in{\bbb N}|\de\le\fracd{\ell (\ell +3)}{2}\}$.
Take $\ell <k$ and assume that
$\ga\in{\bbb P}_{\ell}(P_1,\dots ,P_{\de})$. Since
$$
\{P_{\frac{(\ell +1)\ell}{2}+1},\dots ,P_{\frac{(\ell +2)(\ell +1)}{2}}
\} = L_{\ell +2}\cap (L_1\cup\cdots\cup L_{\ell +1})\sus\ga
$$
one has $\ga =\ga_{\ell -1}\cup L_{\ell +2}$ with
$\ga_{\ell -1}\in{\bbb P}_{\ell -1}(P_1,\dots ,
P_{\frac{(\ell +1)\ell}{2}})$. Continuing this way one finds
$$
\ga =L_{\ell +2}\cup\ga_{\ell -1}=L_{\ell +2}\cup L_{\ell +1}\cup
\ga_{\ell -2}=\cdots =L_{\ell +2}\cup\cdots\cup L_4\cup\ga_1,
$$
where $\ga_j\in{\bbb P}_j(P_1,\dots ,P_{\frac{(j+2)(j+1)}{2}}),
(j=1,\dots ,\ell -1)$. Since $P_1, P_2, P_3$ are not collinear, this is
impossible.

This already proves that ${\bbb P}_{\ell}(P_1,\dots ,P_{\de})=\emptyset$
for $\ell <k$. In particular
${\bbb P}_k(P_1,\dots ,P_{\frac{(k+1)(k+2)}{2}})=\emptyset$. This implies
 $\dim ({\bbb P}_k(P_1,\dots ,P_{\de}))=\fracd{k(k+3)}{2}-\de$. Because
$\de\le\fracd{(d-1)(d-2)}{2}$, $\{ P_1,\dots ,P_{\de}\}\cap L_d=
\emptyset$. So if some element of ${\bbb P}_k(P_1,\dots ,P_{\de})$ would
contain $L_d$ then ${\bbb P}_{k-1}(P_1,\dots ,P_{\de})\ne\emptyset$, a
contradiction.

Hence, ${\bbb P}_k(P_1,\dots ,P_{\de})$ induces a linear system of
dimension $\fracd{k(k+3)}{2}-\de$ on $L_d$. For $P_0$ general on $L_d$
and $\ga\in{\bbb P}_k(P_1,\dots ,P_{\de})$, this implies
$i(\ga, L_d;P_0)\le\fracd{k(k+3)}{2}-\de$, hence
$$
{\bbb P}_k(P_1,\dots ,P_{\de};m)=\emptyset\quad {\rm if\ }m>
\fracd{k(k+3)}{2}-\de .
$$
\vs 1
{\raggedright{\sc Claim}}: There exists a smooth (affine) curve $T$ and
$0\in T$ and a family of plane curves of degree $d$

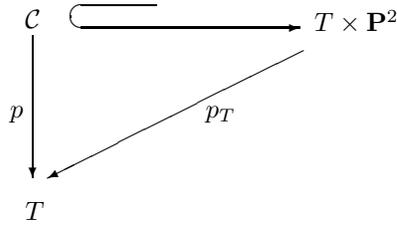
\setlength{\unitlength}{1mm}
\begin{picture}(70,40)(-40,0)
\put(2,6.5){$T$}
\put(0,20){$p$}
\put(26,20){$p_T$}
\put(2,32){${\cal C}$}
\put(40.5,31.5){$T\times{\bbb P}^2$}
\put(3,31){\vector(0,-1){19}}
\put(9.5,33.5){\oval(3,3)[l]}
\put(9.5,32){\vector(1,0){29}}
\put(9.5,35){\line(1,0){10}}
\put(39,29){\vector(-2,-1){34}}
\end{picture}
\newline
with $\de$ sections $S_1,\dots ,S_{\de}:T\to{\cal C}$ satisfying the
following properties:
\begin{enumerate}
\item   $p^{-1}(0)=\Ga_0=L_1\cup\cdots\cup L_d$:
\item   $S_i(0)=P_i$ for $1\le i\le\de$:
\item   for $r\in T-\{ 0\}$, $p^{-1}(r)$ is an irreducible curve,
$S_i(r)$ is an ordinary node for $p^{-1}(r)$ and $p^{-1}(r)$ has no other
 singularities, $P_0$ is a total inflection point on $p^{-1}(r)$.
\end{enumerate}
(For short, we call this a suited family of curves on ${\bbb P}^2$
containing $\Ga_0$ preserving the first $\de$ nodes and the total
inflection point $P_0$.)

Because of semi-continuity reasons it follows that for a general $r\in T$
 the curve $p^{-1}(r)$ satisfies the statement $(*)$. So it is
sufficient to prove the claim.
\vs 2
In order to prove the claim we start as follows. Let
$\pi_1:X_1\to{\bbb P}^2$ be the blowing-up of ${\bbb P}^2$ at $P_0$. Let
$E_1$ be the exceptional divisor and let $L_{d,1}$ be the proper
transform of $L_d$. Let $P^{(1)}=L_{d,1}\cap E_1$. Blow-up $X_1$ at
$P^{(1)}$ obtaining $\pi_2:X_2\to X_1$ with the exceptional divisor $E_2$
 and let $L_{d,2}$ be the proper transform of $L_{d,1}$. Let
$P^{(2)}=L_{d,2}\cap E_2$ and continue until one obtains
$$
\pi:X=X_d\stackrel{\pi_d}{\rightarrow}X_{d-1}
\stackrel{\pi_{d-1}}{\rightarrow}\cdots\stackrel{\pi_2}{\rightarrow}X_1
\stackrel{\pi_1}{\rightarrow}{\bbb P}^2.
$$
Write $L_i$ for $\pi^{-1}(L_i)$ for $1\le i\le d-1$ and let
$$
\Ga'_0=L_1+\dots +L_{d-1}+L_{d,d}.
$$
For $1\le i\le d-1$, let $\mu_i=\pi_{i+1}\circ\cdot\circ\pi_d$ and let
$L$ be a general line on ${\bbb P}^2$. Then
$$
\Ga'_0\in{\bbb P}:=|d\pi^*(L)-\left(\sum_{i=1}^{d-1}\mu^*_i(E_i)\right)
-E_d|
$$

We are going to use a theorem of Tannenbaum \cite[Theorem 2.13]{tannen}.
Since $L_{d,d}.K_X\ge 0$, we are not allowed to take $Y=\Ga'_0$ on $X$ in
 Tannenbaum's Theorem. Therefore we first prove the existence of an
irreducible curve $\Ga'_1$ in ${\bbb P}$ with enough nodes.

{}From Tannenbaum's Theorem it follows that there is a quasi-projective
family ${\bbb P}_d((d-1)(d-2)/2)\sus{\bbb P}_d$ of dimension
$\fracd{d(d+3)}{2}-\fracd{(d-1)(d-2)}{2}$ such that a general element
belongs to a suited family of curves on ${\bbb P}^2$ containing $\Ga_0$
and preserving the first $\fracd{(d-1)(d-2)}{2}$ nodes.

The condition $i(\ga ,L_d;P_0)\ge d$ for $\ga\in{\bbb P}_d((d-1)(d-2)/2)$
 are at most $d$ linear condition. Let
$$
{\bbb P}_d((d-1)(d-2)/2;d)=\{\ga\in{\bbb P}_d((d-1)(d-2)/2)|
i(\ga ,L_d;P_0)\ge d\} .
$$
One has $\Ga_0\in{\bbb P}_d((d-1)(d-2)/2;d)$ and
$$
\dim ({\bbb P}_d((d-1)(d-2)/2;d))\ge\frac{d(d+3)}{2}-\frac{(d-1)(d-2)}{2}
-d=2d-1.
$$
Let $\tilde{\bbb P}$ be an irreducible component of
${\bbb P}_d((d-1)(d-2)/2;d)$ containing $\Ga_0$. Since $\Ga_0$ is smooth
at $P_0$, a general element of $\tilde{\bbb P}$ is smooth at $P_0$. Let
$\Ga_1$ be a general element of $\tilde{\bbb P}$. If $\Ga_1$ is not
irreducible then $i(\Ga_1,L_d;P_0)=d$ implies that $L_d$ is an
irreducible component of $\Ga_1$. Since
$\{ P_1,\dots ,P_{\frac{(d-1)(d-2)}{2}}\}\cap L_d=\emptyset$ also $\Ga_1$
 possesses $\fracd{(d-1)(d-2)}{2}$ nodes none of them belonging to $L_d$.
 This implies $\Ga_1=L_d\cup\Ga_2$, where $\Ga_2$ belongs to a family of
plane curves of degree $d-1$ on ${\bbb P}^2$ containing
$L_2\cup\cdots\cup L_d$ and preserving the $\fracd{(d-1)(d-2)}{2}$ nodes.
 Clearly, if a union of at least two of the lines $L_2,\dots , L_d$
become irreducible in this deformation, some nodes have to disappear.
Since this is not allowed, $\Ga_2$ is the union of $d-1$ lines. But this
would imply $\dim (\tilde{\bbb P})=2d-2$, a contradiction. This proves
that $\Ga_1$ is irreducible.

Moreover $\Ga_1$ belongs to a suited family of curves on ${\bbb P}^2$
containing $\Ga_0$ preserving the first $\fracd{(d-1)(d-2)}{2}$ nodes and
 the total inflection point $P_0$. Because of semi-continuity, we can
assume that $(*)$ holds for the first $\de$ nodes of $\Ga_1$.

Let $\Ga'_1$ be the proper transform of $\Ga_1$ on $X$. Then
$\Ga'_1\in{\bbb P}$ and we can apply Tannenbaum's Theorem to obtain a
suited family of curves on $X$ belonging to ${\bbb P}$ containig $\Ga'_1$
 and preserving the first $\de$ nodes of $\Ga'_1$. Projecting on
${\bbb P}^2$ we obtain a suited family of curves on ${\bbb P}^2$
containing $\Ga_1$, preserving the first $\de$ nodes of $\Ga_1$ and the
total inflection point $P_0$. This completes the proof of the claim.
\vs 2
Let
$$
{\bbb P}_d(d,\de )=\left\{
\begin{array}{ll}
\ga\in{\bbb P}_d: & \ga{\rm\ is\ irreducible;}\\
& \ga{\rm\ has\ a\ total\ inflection\ point\ and}\\
& \ga{\rm\ has\ }\de{\rm\ ordinary\ nodes\ and\ no\ other\ singularities}
\end{array}\right\} .
$$

Then Ran \cite[The irreducibility Theorem (bis)]{ran} proves that
${\bbb P}_d(d;\de)$ is irreducible. This implies:
\vs 1
\begin{thm}
The normalization of a general nodal irreducible plane curve of degree
$d$ with $\de$ nodes and possessing a total inflection point $P$ has in
general Weierstrass gap sequence given by $N^{(1)}_{d,\de}$ at $P$.
\end{thm}

\section{Case: Maximal Weight}
Assume that $\de\le d-2$ and remember the definition for
$N_{d,\de}^{({\rm max})}$ and $N_{d,\de}^{({\rm max2})}$ in the
introduction.

Let $P$ be a total inflection point on the nodal plane curve $\Ga$ of
degree $d$ with $\de$ nodes, let $\al_1 <\cdots <\al_g$ be the
Weierstrass gap sequence of $P$ and let
$N={\bbb N}-\{\al_1,\dots ,\al_g\}$ be the semigroup of non-gaps of $P$.
\vs 1
\begin{lem}
For $1\le i\le g$ one has $\al_i\le\al_i^{({\rm max})}$. Moreover if
$N\ne N_{d,\de}^{({\rm max})}$, then $\al_i\le\al_i^{({\rm max2})}$
for $1\le i\le g$.
\label{lem:max}
\end{lem}
\vs 1
{\it Proof\/}. For $\de\le 2$ see \S 1, so assume that $\de\ge 3$.
Let $\al_{i,j}=(d-i-2)d+j$, $1\le j\le i\le d-2$. They are just the
members of ${\bbb N}-N_d$.

Since $N_d\sus N$, by Lemma \ref{lem:1}, $N$ is the union of $N_d$ and
$\de$ values of $\al_{i,j}$. Moreover, if $\al\in N$ then
$\{\al +d-1,\al + d\}\sus N$. So, if the number of values $\al_{i',j}$
belonging to $N$ with $i'<i$ is less than $\de$, then $\al_{i,j_0}\in N$
for some $1\le j_0\le i$. Each of $N_{d,\de}^{({\rm max})}$ and
$N_{d,\de}^{({\rm max2})}$ does not possess two values
$\al_{i,j_1}\ne\al_{i,j_2}$ for each $i$. Hence, if
$\{\al_{2,1},\al_{2,2}\}\sus N$, then
$$
^{\#}\{\al_{i',j'}\in N|i'<i,j'\ge j\}\ge\ ^{\#}\{
\al_{i',j'}\in N_{d,\de}^{({\rm max2})}|i'<i,j'\ge j\}\quad{\rm for}\
\forall i,j.
$$
So, we have $\al_k\le\al_k^{({\rm max2})}$ for $1\le k\le g$. In
particular, $\al_k\le\al_k^{({\rm max})}$.
But if $\{\al_{2,1},\al_{2,2}\}\not\sus N$, then
$N\in\{ N_{d,\de}^{({\rm max})}, N_{d,\de}^{({\rm max2})}\}$ because
of Corollary \ref{cor:2}.

This completes the proof of the lemma.
\vs 2
\begin{prop}
If $d\ge 2\de+1$, then $N_{d,\de}^{({\rm max})}$ occurs as the semigroup
 of the non-gaps of a total inflection point and if $d\ge 2\de$,
then so does $N_{d,\de}^{({\rm max2})}$.
\label{prop:max}
\end{prop}
\vs 1
{\it Proof\/}. Fix $\de +1$ points $P, P_1,\dots ,P_{\de}$ on an
arbitrary line $L$. For $i=1,\dots ,\de$, take general lines $L_i$ and
$L'_i$ passing through $P_i$. Let $T$ be a general line passing through
$P$ and let $C$ be a curve of degree $d-2\de -1$ which does not pass
through any one of $P, P_1,\dots ,P_{\de}$ and the common point of each
pair of the above curves. Let
\begin{eqnarray*}
C_1 & = & dL\\
C_2 & = & T+C+L_1+L'_1+\cdots +L_{\de}+L'_{\de}.
\end{eqnarray*}
Let ${\bbb P}$ be the pencil generated by $C_1$ and $C_2$. By Bertini's
theorem, a general element $\Ga$ of ${\bbb P}$ is a curve of degree $d$
with $\de$ ordinary nodes at $P_1,\dots ,P_{\de}$ as its only
singularities and $P$ is a total inflection point of $\Ga$ with tangent
line $T$. In particular, if $\Ga$ would not be irreducible then
$\Ga =T+\Ga'$. But then $T$ would be a fixed component of ${\bbb P}$,
which is not true. Hence $\Ga$ is irreducible. Because of Corollary
\ref{cor:2}, the semigroup of nongaps of $P$ is
$N_{d,\de}^{({\rm max})}$.

Next, we prove the latter part. Fix $\de$ points $P_1,\dots ,P_{\de}$ on
an arbitrary line $L$ and a point $P$ not on $L$. For $i=1,\dots ,\de$,
let $L_i$ be the line joining $P$ and $P_i$ and let $L'_i$ be general
lines passing through $P_i$. Let $T$ and $T'$ be general lines passing
through $P$ but not any of $P_i$ and let $C$ be a curve of degree
$d-\de -2$ which does not pass through any one of $P, P_1,\dots ,P_{\de}$
 and the common point of each pair of the above curves. Let
\begin{eqnarray*}
C_1 & = & 2(L_1+\cdots +L_{\de})+(d-2\de )T'\\
C_2 & = & L+T+C+L'_1+\cdots +L'_{\de}.
\end{eqnarray*}
Let ${\bbb P}$ be the pencil generated by $C_1$ and $C_2$. Again, by
Bertini's theorem, a general element $\Ga$ of ${\bbb P}$ is a curve of
degree $d$ with $\de$ ordinary nodes at $P_1,\dots ,P_{\de}$ as its only
singularities and $P$ is a total inflection point of $\Ga$ with tangent
line $T$. Also $\Ga$ is irreducible, by Corollary \ref{cor:2}, the
semigroup of nongaps of $P$ is $N_{d,\de}^{({\rm max2})}$.
\vs 2
{\raggedright {\sc Remark 3.3.}} Define
$N_{d,3}^{({\rm max3})}=N_{d,3}^{({\rm max4})}=N_{d,3}^{(1)}$ and for
$\de >3$ we define inductively $N_{d,\de}^{({\rm max3})}=
N_{d,\de -1}^{({\rm max3})}\cup\{(d-\de -1)d+1\}$ and
$N_{d,\de}^{({\rm max4})}=N_{d,\de -1}^{({\rm max4})}\cup\{(d-\de -1)d
+\de -1\}$. As above one can check that, for $\de\ge 3$ and
$N\not\in\{ N_{d,\de}^{({\rm max})}, N_{d,\de}^{({\rm max2})}\}$ one has
$\al_k\le\al_k^{({\rm max3})}$ for $1\le k\le g$ and, for $\de\ge 5$ and
$N\not\in\{ N_{d,\de}^{({\rm max})}, N_{d,\de}^{({\rm max2})},
N_{d,\de}^{({\rm max3})}\}$ one has $\al_k\le\al_k^{({\rm max4})}$ for
$1\le k\le g$. Moreover $N_{d,\de}^{({\rm max3})}$ (resp.
$N_{d,\de}^{({\rm max4})}$) occurs if and only if exactly $\de -1$ nodes
are on a line $L_0$ and $P\in L_0$ (resp. $P\not\in L_0$). As above one
can also discuss the existence.

If one wants to continue, then one has to start making an analysis of the
 case where the nodes are on a conic. Another direction of further
investigation could be: let $3\le\de'\le\fracd{d}{2}$, what is the
general situation for $N$ if $\de'$ nodes are on a line ?  Probably
reasoning as in \S 2, one obtains an answer.

 \begin{flushleft}
Marc Coppens$^*$\\
Katholieke Industri\"{e}le Hogeschool der Kempen\\
Campus H.~I.~Kempen Kleinhoefstraat 4\\
B 2440 Geel, Belgium\\
\vspace{.1in}
Takao Kato\\
Department of Mathematics,\\
Faculty of Science,\\
Yamaguchi University\\
Yamaguchi, 753 Japan\\
\vspace{.1in}
$^*$This author is related to the University at Leuven\\
(Celestijnenlaan 200B B-3030 Leuven Belgium)\\
as a Research Fellow.
\end{flushleft}
\end{document}